\documentclass{article} 
\usepackage[frenchb,english]{babel} 
\usepackage{amssymb} 
\usepackage{amsmath} 
\usepackage{pb-diagram} 
\usepackage{graphicx} 
\newtheorem{thm}{Theorem}[section]

\newtheorem{pro}[thm]{Proposition} 

\newcommand{\C}{\mathbb{C}} 
\newcommand{\Z}{\mathbb{Z}} 
\newcommand{\R}{\mathbb{R}}

\newcommand{\s}{\mathbb{S}}

\begin{document} 
\title{\bf{A brief note on the spectrum of the basic Dirac operator}} 
\author{{\small{\bf Georges Habib}}\footnote{Max-Planck-Institut f\"ur Mathematik in den Naturwissenschaften, 
Inselstr. 22, 
D-04103 Leipzig, Germany, 
E-mail: Georges.Habib@mis.mpg.de}, {\small\bf{Ken Richardson}}\footnote{Department of Mathematics   
Texas Christian University   
Fort Worth, Texas 76129, USA, E-Mail: k.richardson@tcu.edu}} 
\date{} 
\maketitle 

\begin{abstract} 
In this paper, we prove the invariance of the spectrum of the basic Dirac operator defined on a Riemannian foliation $(M,\mathcal{F})$ with respect to a change of bundle-like metric. We then establish new estimates for its eigenvalues on spin flows 
in terms of the O'Neill tensor and the first eigenvalue of the Dirac operator on $M$. We discuss examples and also define a new version of the basic Laplacian whose spectrum does not depend on the choice of bundle-like metric.   
\end{abstract} 

\section{Introduction and Results} 

\par Let $(M,\mathcal{F})$ be a smooth, closed manifold of dimension $n+1$ endowed with a foliation $\mathcal{F}$ given by an integrable subbundle $L\subset TM$ of rank $p$. We assume throughout the paper that the foliation is Riemannian 
and is endowed with a bundle-like metric \cite{Re}, i.e. the metric $g$ on $M$ induces a holonomy invariant transverse metric on the normal bundle $Q=TM\slash L\simeq L^\perp,$ of rank $q=(n+1)-p;$ this means that the transverse Lie derivative $\mathcal{L}_X g|_Q$ is zero for all leafwise vector fields $X\in\Gamma(L)$. 
This condition is characterized by the existence of a unique metric and torsion-free connection $\nabla$ on $Q$ \cite{Mo}, \cite{Re}, \cite{To}. 
We can then associate to $\nabla$ the transversal curvature data, in particular the transversal Ricci curvature ${\rm Ric}^\nabla$ and transversal scalar curvature ${\rm Scal}^\nabla$. 

\par Many researchers have studied basic forms and the basic Laplacian on Riemannian foliations \cite{Al}, \cite{KT1}, \cite{To}. Basic forms are locally forms on the space of leaves; that is, forms 
$\alpha$ satisfying $X\lrcorner \alpha=X\lrcorner d\alpha=0$ for all $X\in\Gamma(L)$. They are preserved by the exterior derivative and are used to define basic cohomology groups, which are always finite-dimensional in the case of Riemannian foliations. 
The basic Laplacian for a given bundle-like metric is a version of the Laplace operator that preserves the basic forms and 
that is essentially self-adjoint on the $L^2$-closure of the space of basic forms. Its spectrum depends on the choice of the bundle-like metric and provides invariants of that metric. See 
\cite{JuRi}, \cite{LRi1}, \cite{LRi2}, \cite{PaRi}, \cite{Ri1}, \cite{Ri2} for further results. 

\par We now discuss the construction of the basic Dirac operator (see \cite{DGKY}, \cite{GlK}, \cite{PrRi}). 
Let $(M,\mathcal{F})$ be a Riemannian manifold endowed with a Riemannian foliation. 
Let $E\to M$ be a foliated vector bundle (see \cite{KT2}) that is a bundle of $\mathbb{C}\mathrm{l}(Q)$ Clifford modules 
with compatible connection $\nabla^E$.   
The \emph{transversal Dirac operator} $D_{\mathrm{tr}}$ 
is the composition of the maps 
\[ 
\Gamma\left( E\right) \overset{\left( \nabla ^E\right) ^ 
\mathrm{tr}}{\longrightarrow} \Gamma\left( Q^{*}\otimes E\right) \overset{ 
\cong }{\longrightarrow} \Gamma\left( Q\otimes E\right) \overset{{\rm Cliff}}{%
\longrightarrow }\Gamma\left( E\right) , 
\] 
where the last map denotes Clifford multiplication, denoted by $``\cdot",$ and 
the operator $\left( \nabla ^{E}\right) ^{\mathrm{tr}}$ is the 
projection of $\nabla ^{E}.$ 
The transversal Dirac operator fixes the basic sections $\Gamma_b(E)\subset \Gamma(E)$ 
(i.e. $\Gamma_b(E)=\{s\in\Gamma(E):\nabla^{E}_X s=0$ 
for all $X\in \Gamma(L)\}$) but is not symmetric on this subspace. 
By modifying $D_{\mathrm{tr}}$ by a bundle map, we obtain a symmetric and essentially 
self-adjoint operator  $D_{b}$ on $\Gamma_b(E)$. 
To define $D_{b}$, first let 
$H=\sum_{i=1}^p\pi(\nabla^M_{f_i}f_i)$ be the mean curvature of the foliation, where $\pi:TM\rightarrow Q$ denotes the projection and $\{f_i\}_{i=1,\cdots,p}$ is a local orthonormal frame of $L$. Let $\kappa=H^\flat$ be the corresponding $1$-form, so that $H=\kappa^\sharp$. Let $\kappa_b:=P_b \kappa$ be the $L^2$-orthogonal 
projection of $\kappa$ onto the space of basic forms (see \cite{Al}, \cite{PaRi}), and let 
$\kappa_b^\sharp$ be the corresponding vector field. 
We now define 
\begin{equation*} 
D_{b}s:=\frac{1}{2}(D_{{\rm tr}}+D_{{\rm tr}}^*)s=\sum_{i=1}^q e_i\cdot \nabla^E_{e_i}s-\frac{1}{2}\kappa_b^{\sharp}\cdot s , 
\label{eq:bas} 
\end{equation*} 
where $\{e_i\}_{i=1,\cdots,q}$ is a local orthonormal frame of $Q$. 
A direct computation shows that $D_{b}$ preserves the basic sections, 
is transversally elliptic, and thus has discrete spectrum \cite{EK}. 

\par In this paper, we study the invariance of the spectrum of the basic Dirac operator with respect to a change of bundle-like metric; that means when one modifies the metric on $M$ in any way that leaves the transverse metric on the normal bundle intact (this includes modifying the subbundle $L^{\perp}\subset TM$, as one must do in order to make the mean curvature basic, for example). We prove 
\begin{thm}\label{inv} Let $(M,\mathcal{F})$ be a compact Riemannian manifold endowed with a Riemannian foliation and basic Clifford bundle $E\to M$. The spectrum of the basic Dirac operator is the same for every possible choice of bundle-like metric that is associated to the transverse metric on the quotient bundle $Q$. 
\end{thm} 

\par We emphasize that the basic Dirac operator $D_b$ depends on the choice of bundle-like metric, not merely on the Clifford structure and Riemannian foliation structure, since both projections $T^*M\to Q^*$ and $P_b$ as well as $\kappa_b$ depend on the leafwise metric. We clarify this result on the $2$-dimensional torus (see Example $2$ at the end of Section \ref{pro}). In Section \ref{basLap}, we show that the square of the standard 
basic Dirac operator on basic forms is related to but not the same as the basic Laplacian, 
whose spectrum does depend on the choice of compatible bundle-like metric. 

\par In \cite{Hab}, the first author studied the spectrum of the basic Dirac operator on spin foliations. That is, a compact Riemannian manifold $(M,\mathcal{F})$ endowed with a Riemannian foliation $\mathcal{F}$ of codimension $q$ such that the normal bundle carries a spin structure (as a foliated vector bundle). In that paper, the author obtained a Friedrich-type estimate \cite{Fr} for the eigenvalues of the basic Dirac operator for bundle-like metrics with basic-harmonic mean curvature. Since by Theorem \ref{inv} the spectrum does not change, 
we deduce that for any bundle-like metric the eigenvalues of the basic Dirac operator satisfy   
\begin{equation} 
\lambda^2\geq \frac{q}{4(q-1)}{\rm inf}_M({\rm Scal}^\nabla). 
\label{eq:esti} 
\end{equation} 
This estimate is of interest only for positive transversal scalar curvature. Moreover, we show 
\begin{pro} \label{pro:esti} Let $(M,\mathcal{F})$ be a compact Riemannian manifold endowed with a spin foliation with basic mean curvature $\kappa.$ Then, we have the estimate 
\begin{equation*} 
\lambda^2\geq \frac{q}{4(q-1)}{\rm inf}_M({\rm Scal}_M-{\rm Scal}_L+|A|_Q^2+|T|_L^2). 
\label{eq:estima} 
\end{equation*}   
If $\mathcal{F}$ is a Riemannian flow (i.e. $p=1$), then 
\begin{equation} 
\lambda^2\geq \frac{q}{4(q-1)}{\rm inf}_M({\rm Scal}_M+|A|_Q^2+|\kappa|^2). 
\label{eq:estmflot} 
\end{equation} 
If the limiting case is attained, the foliation is minimal and we have a transversal Killing spinor. 
\end{pro} 
Here $A$ and $T$ denote the O'Neill tensors \cite{B,O} of the foliation (see Section \ref{pro} for the definitions). We point out that Estimate \eqref{eq:estmflot} may improve Inequality \eqref{eq:esti} (see Example $1$). 
On the other hand, using the Min-Max principle we give a new lower bound for the eigenvalues of the basic Dirac operator on a Riemannian flow $(M,\mathcal{F})$ in terms of the first eigenvalue of the Dirac operator of $M$ and the O'Neill tensor of the flow. We prove 
   
\begin{pro} \label{pro:esti1} Let $(M^{n+1},\mathcal{F})$ be a compact spin manifold endowed with a Riemannian flow $\mathcal{F}$ given by a unit vector field $\xi$. Then we have the estimate   
\begin{equation} 
\lambda^2\geq \frac{1}{2}\lambda^2(D_M)-\frac{n}{16}{\rm sup}_M(|A|_Q^2), 
\label{eq:5} 
\end{equation}   
where $\lambda(D_M)$ is the first eigenvalue of the Dirac operator of $M$. 
\end{pro} 
   
\section{Proof of Theorem \ref{inv}} 
\par Suppose that  $D_b$ is any basic Dirac operator associated to a basic Clifford bundle over a Riemannian foliation with bundle-like metric. 
Consider a change of the bundle-like metric as mentioned in the introduction. When one does this, there is no effect on the compatibility of the new transversal connection, since the bundles and bundle metrics are unchanged. Furthermore, the set of all basic sections of the Clifford bundle is unchanged. Let $D_b'$ denotes the new basic Dirac operator, $\mathrm{dvol}'$ the new volume form, and $\kappa_b'$ the basic component of mean curvature with respect to the new bundle-like metric. We see that 
$\mathrm{dvol}'= h\, \mathrm{dvol}$, for some positive smooth function $h$. Then the new $L^2$-metric on basic sections $s$, $t$ is 
\[ 
<s,t>' = \int_M (s,t) h \,\mathrm{dvol}, 
\] 
where $(\,\, , \,\,)$ is the pointwise inner product on $\Gamma(E)$. As in \cite{PaRi}, we have 
\[ 
<s,t>' = \int_M (s,t) P_b h \,\mathrm{dvol}. 
\] 
Setting $\alpha=P_b h$, we see that 
\begin{eqnarray*} 
<\alpha^{-1/2}D_b(\alpha^{1/2}s),t>' &=& \int_M ( \alpha^{-1/2}D_b(\alpha^{1/2}s) , t ) \alpha \,\mathrm{dvol}\\ 
&=& \int_M ( D_b(\alpha^{1/2}s) , \alpha^{1/2}t ) \,\mathrm{dvol}\\ 
&=& \int_M ( \alpha^{1/2}s , D_b(\alpha^{1/2}t) ) \,\mathrm{dvol}, 
\end{eqnarray*} 
where the last equality is a consequence of the fact that $D_b$ is self-adjoint with respect to the old $L^2$-metric. Hence, 
\begin{eqnarray*} 
<\alpha^{-1/2}D_b(\alpha^{1/2}s),t>' &=& \int_M ( s , \alpha^{-1/2}D_b(\alpha^{1/2}t) ) \alpha \,\mathrm{dvol} \\ 
&=& < s ,\alpha^{-1/2}D_b(\alpha^{1/2}t) >' , 
\end{eqnarray*} 
so that $\alpha^{-1/2}D_b \alpha^{1/2}= D_b + \frac{1}{2\alpha}(d\alpha)\cdot$ 
is self-adjoint with respect to the new $L^2$-metric. Now, we compute 
\begin{eqnarray*} 
D_b'-\alpha^{-1/2}D_b \alpha^{1/2}&=& D_b'-D_b - \frac{1}{2\alpha}(d\alpha)\cdot \\ 
&=& D_{\mathrm{tr}}-\frac{1}{2}\kappa_b'^{\sharp}\cdot-D_{\mathrm{tr}}+\frac{1}{2}\kappa_b^{\sharp}\cdot- \frac{1}{2\alpha}(d\alpha)\cdot \\ 
&=&-\frac{1}{2}\kappa_b'^{\sharp}\cdot+\frac{1}{2}\kappa_b^{\sharp}\cdot- \frac{1}{2\alpha}(d\alpha)\cdot. 
\end{eqnarray*} 
By the symmetry of $D_b'-\alpha^{-1/2}D_b \alpha^{1/2}$ with respect to the new $L^2$-metric, 
the operator $\kappa_b'^{\sharp}\cdot-\kappa_b^{\sharp}\cdot- \frac{1}{2\alpha}(d\alpha)\cdot$ must be 
a Hermitian symmetric endomorphism of $\Gamma(E)$. Since Clifford multiplication by cotangent 
vectors is skew-Hermitian, this bundle endomorphism must be zero. 
Since $D_b'$ is just a conjugate of $D_b$, all the spectrum is unchanged. 
\hfill$\square$\\ 

\par{\bf Remark.} Dom\'inguez showed in \cite{Do} that every Riemannian foliation admits a bundle-like metric with basic mean curvature 
(so that  $\kappa=\kappa_b$). Further, Mason showed in \cite{Ma} that the bundle-like metric may be chosen so that the mean curvature is basic-harmonic ($\Delta_b\kappa=0$). Thus, we may assume that the mean curvature form has this property when calculating the eigenvalues of the basic Dirac operator. 

\section{Proof of Propositions \ref{pro:esti} and \ref{pro:esti1}} \label{pro} 
We recall that the geometry of the normal bundle of a Riemannian foliation is characterized by the existence of the O'Neill tensors, given for all $X,Y\in \Gamma(TM)$ by 
\begin{equation*} 
A_X Y=\pi^\perp(\nabla^M_{\pi(X)}\pi(Y))+\pi(\nabla^M_{\pi(X)}\pi^\perp(Y)), 
\end{equation*} 
and 
\begin{equation*} 
T_X Y=\pi^\perp(\nabla^M_{\pi^\perp(X)}\pi(Y))+\pi(\nabla^M_{\pi^\perp(X)}\pi^\perp(Y)), 
\end{equation*} 
where $\pi^\perp:TM\rightarrow L$ denotes the projection. Indeed, for all $Z\in \Gamma(Q)$ the tensor $A_Z$ is a skew-symmetric map interchanging $L$ and $Q$ and we have that $A_Z W=\frac{1}{2}\pi^\perp[Z,W]$ for all $W\in \Gamma(Q)$. Moreover, the scalar curvatures of $M$, $Q$ and $L$ are related by the following relation \cite[p. 244]{B} 
\begin{equation} 
{\rm Scal}_M={\rm Scal}^\nabla+{\rm Scal}_L-|A|_Q^2-|T|_L^2-|\kappa|^2+2{\rm div}_Q\kappa 
\label{eq:scal} 
\end{equation} 
where $|\cdot|_Q$ (resp. $|\cdot|_L$) is the norm evaluated on a local frame of $\Gamma(Q)$ (resp. $\Gamma(L)$).\\   

\par For the proof of Proposition \ref{pro:esti}, we assume that the mean curvature is basic and use the transversal Schr\"odinger-Lichnerowicz formula for $D_b^2$ \cite[Eq. 1.6]{GlK} to write for any basic spinor $\psi$ 
\begin{equation} 
D_b^2\psi=\nabla^*\nabla\psi+\frac{1}{4}({\rm Scal}^\nabla+|\kappa|^2)\psi-\frac{1}{2}\delta_b\kappa\cdot\psi , 
\label{eq:schlich} 
\end{equation} 
where $\delta_b$ is the adjoint of $d$ on the space of basic forms \cite{PaRi}. Plugging Equation \eqref{eq:scal} in \eqref{eq:schlich} finishes the proof of the proposition after integrating over $M$ the scalar product of \eqref{eq:schlich} by $\psi$ and using the fact that 
$\mathrm{div}_Q \kappa = -\delta_b \kappa + |\kappa |^2$ 
if $\kappa$ is basic and 
that $|\nabla\psi|^2\geq\frac{1}{q}\lambda^2|\psi|^2$. 
\hfill$\square$\\ 

\par{\bf Remark.} In \cite{Hab}, we consider a Riemannian flow $\mathcal{F}$ and modify the metric of $M$ in the direction of the flow, i.e. we consider the family of metrics $g_t=t^2g_\xi\oplus g_Q$. The spectrum of the Dirac operator of $(M,g_t)$ tends to the spectrum of the basic Dirac operator $D_b$ when $t$ tends to zero. Hence by Friedrich's estimate, we have for the Dirac operator of $(M,g_t)$ that $\lambda_t^2\geq \frac{q+1}{4q}\inf_M{\rm Scal}_M^t$. A direct computation for the scalar curvature leads to ${\rm Scal}_M^t={\rm Scal}_M+(1-t^2)|A|_Q^2$. Thus when $t$ tends to zero, we obtain the estimate for the eigenvalues of the basic Dirac operator $\lambda^2\geq \frac{q+1}{4q}\inf_M{\rm (Scal}_M+|A|_Q^2)$ which is not as strong as \eqref{eq:estmflot} since $\frac{q}{(q-1)}\geq\frac{q+1}{q}.$\\ 
   
\par The proof of Proposition \ref{pro:esti1} is based on the use of the Min-Max principle by computing the Rayleigh quotient  $\frac{||D_M\psi||^2}{||\psi||^2}$ associated with a basic spinor $\psi,$ where $||\cdot||$ 
is the $L^2$-norm on the spinor bundle $\Sigma M$ induced by the Hermitian metric. For this, we recall that the spinorial Levi-Civita connections $\nabla^M$ and $\nabla$ are related for any spinor field $\psi\in \Gamma(\Sigma M)$ by \cite[Eq. 2.4.7]{Hab} 
\begin{equation*} 
\left\{\begin{array}{ll}\nabla_\xi^M\psi&= \nabla_\xi\psi+\frac{1}{4}\sum_{i=1}^n e_i\cdot A_{e_i}\xi\cdot\psi+\frac{1}{2}\xi\cdot\kappa\cdot\psi\\ 
 &\\ 
\nabla_Z^M\psi&= \nabla_Z\psi+\frac{1}{2}\xi\cdot A_Z\xi\cdot\psi,\end{array}\right. 
\end{equation*} 
where $Z\in \Gamma(Q)$ and $\{e_i\}_{i=1,\cdots,n}$ is an orthonormal frame of $\Gamma(Q)$. 
We easily deduce the relation between the basic Dirac operator and the Dirac operator of $M$: 
$$D_M\psi=D_b\psi-\frac{1}{4}\sum_{i=1}^n\xi\cdot e_i\cdot A_{e_i}\xi\cdot\psi,$$ 
for each basic spinor field $\psi$. Using the pointwise inequality $|\sum_{i=1}^n\xi\cdot e_i\cdot A_{e_i}\xi\cdot\psi|\leq \sum_{i=1}^n |A_{e_i}\xi| |\psi|$ and the fact that $|a+b|^2\leq 2a^2+2b^2,$ we deduce that for any basic eigenspinor $\psi$ with eigenvalue $\lambda$, we have 
\begin{eqnarray*} 
|D_M\psi|^2&\leq &2\lambda^2|\psi|^2+2(\frac{1}{4})^2(\sum_{i=1}^n|A_{e_i}\xi|)^2|\psi|^2\nonumber\\ 
&\leq & 2\lambda^2|\psi|^2+\frac{n}{8}|A|_Q^2|\psi|^2, 
\end{eqnarray*} 
where the last inequality is a consequence of the Cauchy-Schwarz inequality. Thus we get that $\frac{||D_M\psi||^2}{||\psi||^2}\leq 2\lambda^2+\frac{n}{8}{\rm sup}_M(|A|_Q^2)$, 
and by the Min-Max principle the proof of the proposition follows. 
\hfill$\square$\\ 
   
\par Now we will consider particular Riemannian flows on the $3$-dimensional sphere and the $2$-dimensional torus and will estimate the eigenvalues of the basic Dirac operator. \\ 

\noindent{\bf Example 1.} Let $\s^3\subset\C^2$ be the round $3$-sphere equipped with its standard metric of sectional curvature equal to $1$ and its standard spin structure. For $r>0$, we define the $1$-parameter group $\gamma_t^r$ of isometries by \cite{Heb} 
$$\gamma_t^r(z,w)=(e^{irt}z,e^{it}w)$$ 
for $t\in\R$ and $(z,w)\in \s^3$. This group generates the Killing vector field $Z^r\in T\s^3$ defined by $Z^r_{(z,w)}=(irz,iw)$ and thus the vector field $\xi^r_{(z,w)}=Z^r_{(z,w)}/|Z^r_{(z,w)}|$ defines a Riemannian flow $\mathcal{F}_r$ on $\s^3$. 
Let us consider on $T_{(z,w)}\s^3$ the two vector fields $X^r_{(z,w)}=(|w|^2z,-|z|^2w)$ and $Y^r_{(z,w)}=(i|w|^2z,-ir|z|^2w)$ for $(z,w)\in \s^3$. The frame   
$$\{\xi^r=Z^r/|Z^r|,\,\,\, e_1^r=X^r/|X^r|,\,\,\, e_2^r=Y^r/|Y^r|\},$$ 
is orthonormal of the tangent space of $\s^3$. 
Hence by a direct computation with respect to this frame, on gets that the transversal scalar curvature is equal to ${\rm Scal}^{\nabla}_{(z,w)}=2+\frac{6r^2}{r^2|z|^2+|w|^2}$ and the mean curvature $\kappa^r$ of the flow is a basic $1$-form and equal to $\frac{(1-r^2)|z||w|}{r^2|z|^2+|w|^2}e^r_1$. The O'Neill tensor of the flow is $A_{e_1^r}e_2^r=\frac{r}{r^2|z|^2+|w|^2}\xi^r$. 

\par On the one hand, Inequality \eqref{eq:esti} allows to deduce an estimate for the eigenvalues of the basic Dirac operator 
$$ 
\lambda^2>\frac{1}{2}{\rm inf}_{\s^3}(2+6\frac{r^2}{r^2|z|^2+|w|^2})=\left\{\begin{array}{ll} 
1+3r^2 &  \textrm {if $r<1$}\\\\ 
4 & \textrm {if $r>1$}. 
\end{array}\right.$$   
On the other hand, by Inequality \eqref{eq:estmflot} we get that 
$$\lambda^2> \frac{1}{2}{\rm inf}_{\s^3} (6+\frac{2r^2+(1-r^2)^2|z|^2|w|^2}{(r^2|z|^2+|w|^2)^2})=\left\{\begin{array}{ll} 
r^2+3 &  \textrm {if $r<1$}\\\\ 
\frac{1}{r^2}+3 & \textrm {if $r>1$}. 
\end{array}\right.$$   
We see that for $r<1$ the last estimate is better than the first one. The lower bound in Inequality \eqref{eq:5} is weaker than those above since we get that $\lambda^2> \left\{\begin{array}{ll} 
\frac{9}{8}-\frac{1}{4r^2} &  \textrm {if $r<1$}\\ 
\frac{9}{8}-\frac{r^2}{4}& \textrm {if $r>1$}. 
\end{array}\right.$\\\\ 

\noindent{\bf Example 2.} Consider the manifold $(\s^1\times \s^1, f^2d\theta^2\oplus dt^2)$ where $f(\theta,t)$ is a function 
on $\s^1\times \s^1=\mathbb{R}^2\slash(2\pi\mathbb{Z})^2$. 
The frame $\xi=\frac{1}{f}\partial_\theta, e_1=\partial_t$ is orthonormal with respect to the metric and we have 
$[\xi,e_1]=\frac{1}{f}\frac{\partial f}{\partial t}\xi$. The vector field $\xi$ defines a Riemannian flow on $M$ with vanishing O'Neill tensor and mean curvature is $\kappa^{\sharp}=-\frac{1}{f}\frac{\partial f}{\partial t}e_1$. 
The basic component of the mean curvature satisfies $\kappa_b^{\sharp}=-\frac{\dot{g}(t)}{g(t)}e_1$, 
where $g(t)=\frac{1}{2\pi}\int_0^{2\pi}f(\theta,t)\, d\theta$. 

\par In order to compute the spectrum of the basic Dirac operator, we will find solutions of the differential equation $D_b\psi=\lambda\psi$, with $\lambda$ is an eigenvalue of $D_b$. Clifford multiplication by $e_1$ 
is multiplication by $i$. Hence the equation $D_b\psi=\lambda\psi$ can be written as 
$$\frac{\partial}{\partial \theta}\psi=0 \quad\text{and}\quad i \frac{\partial}{\partial t}\psi+i\frac{\dot{g}(t)}{2g(t)}\psi=\lambda \psi.$$ 
The solution of this equation is $\psi(t)=cg^{-\frac{1}{2}}(t) e^{-i\lambda t}$, where $c$ is a complex constant. For the trivial spin structure, the spinor $\psi$ is $2\pi$-periodic if $\lambda=k\in \Z.$ This verifies that the spectrum does 
not depend on the choice of the function $f$, as Theorem \ref{inv} shows. 

\par This example may be easily modified to provide a similar example on any suspension of an isometry $\varphi:M\rightarrow M$, yielding a one-dimensional foliation of an $M$-bundle over a circle. That is, the manifold is  $\mathbb{R}\times_\varphi M$, 
 the quotient of $\mathbb{R}\times M$ by the equivalence relation 
 $(\,\theta,x\,)\sim (\,\theta-2\pi,\varphi(x) \,)$.

\section{Remarks about the basic Laplacian} \label{basLap} 
As mentioned in the introduction, it is well-known that the eigenvalues of the basic Laplacian depend on the choice of bundle-like metric; for example, in \cite[Corollary 3.8]{Ri2}, it is shown that the spectrum of the basic Laplacian on functions determines the $L^2$-norm of the mean curvature on a transversally oriented foliation of codimension one. If the foliation were taut, then a bundle-like metric could be chosen so that the mean curvature is identically zero, and other metrics could be chosen where the mean curvature is nonzero. This is one reason why the invariance of the spectrum of the basic Dirac operator is such 
a surprise. 

\par If one considers the bundle $\bigwedge^*Q^*$ as a $\mathbb{C}\mathrm{l}(Q)$ bundle with Clifford multiplication 
\[ 
\alpha^\sharp\cdot = \alpha\wedge - \alpha \lrcorner 
\] 
for (local) basic one-forms $\alpha$, then the corresponding basic Dirac operator is 
$$D_b=d+\delta_b-\kappa_b\lrcorner-\frac{1}{2}\kappa_b^\sharp\cdot.$$ 
The square of this operator is certainly not the basic Laplacian $\Delta_b$, but 
$\Delta_b$ and $D_b^2$ do have the same principal symbol. 
Observe that since  $d\kappa_b=0$ (see \cite{Al}, \cite{PaRi}), 
the following operator is a differential on the space of basic forms: 
\[ 
\widetilde{d}:=d-\frac{1}{2}\kappa_b\wedge, 
\] 
and its $L^2$-adjoint restricted to the space of basic forms is 
\[ 
\widetilde{\delta}_b=\delta_b-\frac{1}{2}\kappa_b\lrcorner\, . 
\] 
Then 
\[ 
\widetilde{d}+\widetilde{\delta}_b=d+\delta_b-\frac{1}{2}\kappa_b\lrcorner-\frac{1}{2}\kappa_b\wedge=D_b. 
\] 
Thus, with this new twisted differential on basic forms, the analogy with ordinary manifolds carries over, 
and the resulting basic $D_b^2$-harmonic forms would play the same role 
in Hodge theory with the twisted basic $\widetilde{d}$-cohomology. Further, unlike the ordinary and well-studied basic 
Laplacian, the eigenvalues of  $D_b^2$ are invariants of the Riemannian foliation 
structure alone and independent of the choice of compatible bundle-like metric.

\end{document}